\documentclass[centertags,11pt]{article}

\usepackage{amsmath}
\usepackage{amsthm}
\usepackage{amscd}
\usepackage{amssymb}
\usepackage{verbatim}

\title{Isometries, rigidity and universal covers}
\author{Benson Farb and Shmuel Weinberger \thanks{Both authors are
supported in part by the NSF.}}

\numberwithin{equation}{section}

\theoremstyle{plain}
\newtheorem{theorem}{Theorem}[section]
\newtheorem{proposition}[theorem]{Proposition}

\newtheorem{corollary}[theorem]{Corollary}
\newtheorem{claim}[theorem]{Claim}
\newtheorem{conjecture}[theorem]{Conjecture}
\newtheorem{question}[theorem]{Question}

\theoremstyle{definition}

\newtheorem{remark}[theorem]{Remark}
\newtheorem{example}[theorem]{Example}

\def\proof{{\bf {\medskip}{\noindent}Proof. }}

\def\endproof{$\diamond$ \bigskip}

\def\title{\em}

\newcommand\R{\mbox{\bf R}}

\newcommand\Z{\mbox{\bf Z}}
\newcommand\Q{\mbox{\bf Q}}

\DeclareMathOperator{\Out}{Out}

\DeclareMathOperator{\qcd}{cd_{\Q}}

\DeclareMathOperator{\rank}{rank}

\DeclareMathOperator{\U}{U}
\DeclareMathOperator{\GL}{GL}

\DeclareMathOperator{\Isom}{Isom}

\DeclareMathOperator{\Aut}{Aut}

\DeclareMathOperator\nbhd{Nbhd}

\DeclareMathOperator\cdq{cd_{\bf Q}}
\DeclareMathOperator\Comm{Comm}

\DeclareMathOperator\frames{Fr}
\DeclareMathOperator\sol{sol}
\DeclareMathOperator\sem{ss}
\DeclareMathOperator\vol{Vol}

\renewcommand\to{\longrightarrow}

\begin{document}
\maketitle

\section{Introduction}
The goal of this paper is to describe all closed, aspherical Riemannian
manifolds $M$ whose universal covers $\widetilde{M}$ have have a
nontrivial amount of symmetry.  By this we mean that $\Isom(\widetilde{M})$
is not discrete.  By the well-known theorem of Myers-Steenrod \cite{MS},
this condition is equivalent to
$[\Isom(\widetilde{M}):\pi_1(M)]=\infty$.  Also note that if any cover
of $M$ has a nondiscrete isometry group, then so does its universal
cover $\widetilde{M}$.  

Our description of such $M$ is given in 
Theorem \ref{theorem:main} below.  The proof of this theorem 
uses methods from Lie theory, 
harmonic maps, large-scale geometry, and 
the homological theory of transformation groups.  

The condition that $\widetilde{M}$ have nondiscrete isometry group 
appears in a wide variety of problems in geometry.  Since Theorem
\ref{theorem:main} provides a taxonomy of such $M$, it can be used 
to reduce many general problems to verifications of specific examples.
Actually, it is not always Theorem \ref{theorem:main} which is applied
directly, but the main subresults from its proof.  After
explaining in \S\ref{explain} the statement of Theorem \ref{theorem:main}, we 
give in \S\ref{section:intro:cors} a number of such applications.  These
range from
new characterizations of locally symmetric manifolds, to the
classification of contractible manifolds covering both compact and
finite volume manifolds, to a new proof of the Nadel-Frankel Theorem 
in complex geometry.   

\subsection{Statement of the general theorem}
\label{explain}

The basic examples of closed,
aspherical, Riemannian manifolds whose universal covers have 
nondiscrete isometry groups are the {\em locally homogeneous
(Riemannian) manifolds} $M$, 
i.e. those $M$ whose universal cover admits a transitive Lie group
action whose isotropy subgroups are maximal compact.  Of
course one might also take a product of such a manifold with an
arbitrary manifold.  To find nonhomogeneous examples which are not
products, one can do the following construction.

\begin{example}
Let $F\rightarrow M\rightarrow B$ be any Riemannian fiber bundle with
the induced path metric on $F$ locally homogeneous.  Let $f:B\rightarrow
\R^+$ be any smooth function.  Now at each point of $M$ lying over $b$,
rescale the metric in the tangent space $TM_b=TF_b\oplus TB_b$ by 
rescaling $TF_b$ by $f(b)$.  Almost any $f$ gives a metric on
$M$ with $\dim(\Isom(\widetilde{M}))>0$ but with $\widetilde{M}$
not homogeneous, indeed with each $\Isom(\widetilde{M})$-orbit a
fiber.  This construction can be further extended by scaling
fibers using any smooth map from $B$ to the moduli space of
locally homogeneous metrics on $F$; this moduli space is large for
example when $F$ is an $n$-dimensional torus.
\end{example}

Hence we see that there are many closed, aspherical, Riemannian
manifolds whose universal covers admit a nontransitive action of a
positive-dimensional Lie group.  
The following general result says that the examples described above 
exhaust all the possibilities for such manifolds.  

Before stating the general result, we need some terminology.  
A {\em Riemannian orbifold} $B$ is a smooth orbifold where
the local charts are modelled on quotients $V/G$, where $G$ is a finite
group and $V$ is a linear
$G$-representation endowed with some $G$-invariant Riemannian metric.  
The orbifold $B$ is {\em good} if it is the quotient of $V$
by a properly discontinuous group action.  

A {\em Riemannian orbibundle} is a smooth map
$M\to B$ from a Riemannian manifold to a Riemannian orbifold locally
modelled on the quotient map $p:V\times_GF\to V/G$, where $F$ is a fixed
smooth manifold with smooth $G$-action, and where $V\times F$ has a
$G$-invariant Riemannian metric such that projection to $V$ is an
orthogonal projection on each tangent space.  Note that in this
definition, the induced metric on the fibers of a Riemannian orbibundle
may vary, and so a Riemannian orbibundle is not a fiber bundle structure
in the Riemannian category.

\begin{theorem}
\label{theorem:main}
Let $M$ be a closed, aspherical Riemannian manifold.  Then either
$\Isom(\widetilde{M})$ is discrete, or $M$ is isometric to an orbibundle 
\begin{equation}
\label{eq:general:orbi}
F\to M\to
B
\end{equation} 
where:
\begin{itemize}
\item $B$ is a good Riemannian orbifold.
\item Each fiber $F$, endowed with the induced metric, is isometric to a
closed, aspherical, locally homogeneous Riemannian $n$-manifold, $n>0$
\footnote{Recall that a manifold $F$ is {\em locally homogeneous} if its
universal cover is isometric to $G/K$, where $G$ is a Lie group, $K$ is
a maximal compact subgroup, and $G/K$ is endowed with a left
$G$-invariant, $K$ bi-invariant metric.}.
\end{itemize}

\end{theorem}

Note that $B$ is allowed to be a single point.  

One might hope that the Riemannian orbifold $B$ in the conclusion of Theorem
\ref{theorem:main} could be taken to be a Riemannian manifold, at least
after passing to a finite cover of $M$.  This is not the case,
however.  In \S\ref{section:example:nocover} we construct a Riemmanian 
manifold $M$ with the property that $M$ is a Riemannian
orbibundle, fibering over a singular orbifold, 
but such that no finite cover of $M$ fibers over a manifold; further, 
$\Isom(\widetilde{M})$ is not discrete.  This
seems to be the first known example of an aspherical manifold with 
a singular fibration that 
remains singular in every finite cover. In constructing $M$ we 
produce a group $\Gamma$ acting properly discontinuously and cocompactly
by diffeomorphisms on $\R^n$, but which is not virtually torsion-free. 

\subsection{Applications}
\label{section:intro:cors}

We now explain how to apply Theorem \ref{theorem:main} and its proof 
to a variety 
of problems in geometry.  The proofs of 
these results will be given in \S\ref{section:consequences} below.  

\bigskip
\noindent
{\bf Characterizations of locally symmetric manifolds. }We begin
with a characterization of locally symmetric manifolds among all
closed Riemannian manifolds.  The theme is that such manifolds are
characterized by some simple properties of their fundamental group,
together with the property that their universal covers have 
nontrivial symmetry (i.e.\ have nondiscerete isometry group).  
We say that a smooth manifold $M$ is {\em
smoothly irreducible} if $M$ is not smoothly covered by a 
nontrivial finite product of smooth manifolds.

\begin{theorem}
\label{theorem:locsym1}
Let $M$ be any
closed Riemannian $n$-manifold, $n>1$.  Then the following are equivalent:
\begin{enumerate}
\item $M$ is aspherical, smoothly irreducible, 
$\pi_1(M)$ has no nontrivial, normal abelian 
subgroup, and $\Isom(\widetilde{M})$ is not discrete.
\item $M$ is isometric to an irreducible, locally-symmetric Riemannian
manifold of nonpositive sectional curvature.
\end{enumerate}
\end{theorem}

The idea here is to apply Theorem \ref{theorem:main}, or
more precisely the main results in its proof, and then 
to show that if the base $B$ were positive dimensional, the manifold 
$M$ would not be smoothly irreducible; see
\S\ref{section:proof:locsym} below.

\medskip
\noindent
{\bf Remark. }The proof of Theorem \ref{theorem:locsym1} gives more: 
the condition that $M$ is smoothly irreducible 
can be replaced by the weaker condition that $M$ is not 
Riemannian covered by a nontrivial Riemannian warped product; see 
\S\ref{section:proof:locsym}.
\bigskip

When $M$ has nonpositive curvature, the Cartan-Hadamard Theorem gives
that $M$ is aspherical.  For nonpositively curved metrics on $M$, 
Theorem \ref{theorem:locsym1} was proved by Eberlein in
\cite{Eb1,Eb2} 
\footnote{Eberlein's results are proved not just for lattices but more
generally for groups satisfying the so-called {\em
duality condition} (see \cite{Eb1,Eb2}), a condition on the limit set of
the group acting on the visual boundary.}.  While the differential geometry
and dynamics related to nonpositive curvature are central to 
Eberlein's work, for the most part they do not, by necessity, play a
role in this paper.

Recall that the Mostow Rigidity Theorem states that a closed, aspherical
manifold of dimension at least three admits at most one 
irreducible, nonpositively curved, locally symmetric metric up to homotheties
of its local direct factors.  For such locally symmetric manifolds $M$, 
Theorem \ref{theorem:locsym1} has the following immediate consequence:

\medskip
{\it Up to homotheties of its local direct factors, the 
locally symmetric metric on $M$ is the unique Riemannian 
metric with $\Isom(\widetilde{M})$ not discrete.}

\medskip

Uniqueness within the set of nonpositively curved Riemannian metrics on
$M$ follows from \cite{Eb1,Eb2}.  This statement also generalizes
the characterization in \cite{FW} of the locally symmetric metric on an 
arithmetic manifold.

Combined with basic facts about word-hyperbolic groups, Theorem
\ref{theorem:locsym1} provides the following characterization of closed,
negatively curved, locally symmetric manifolds.

\begin{corollary}
\label{corollary:locsym3}
Let $M$ be any closed Riemannian $n$-manifold, $n>1$. Then the following are
equivalent:
\begin{enumerate}
\item $M$ is aspherical, $\pi_1(M)$ is word-hyperbolic, and 
$\Isom(\widetilde{M})$ is not discrete.
\item $M$ is isometric to a negatively curved, 
locally symmetric Riemannian manifold.
\end{enumerate}
\end{corollary}

Theorem \ref{theorem:locsym1} can also be combined with Margulis's
Normal Subgroup Theorem to give a simple characterization in the higher 
rank case.  We say that a
group $\Gamma$ is {\em almost simple} if every normal subgroup of
$\Gamma$ is finite or has finite index in $\Gamma$.

\begin{corollary}
\label{corollary:locsym2}
Let $M$ be any closed Riemannian manifold. Then the following are
equivalent:
\begin{enumerate}
\item $M$ is aspherical, $\pi_1(M)$ is almost simple, and 
$\Isom(\widetilde{M})$ is not discrete.
\item $M$ is isometric to a nonpositively curved, irreducible, 
locally symmetric Riemannian manifold of (real) rank at least $2$.
\end{enumerate}
\end{corollary}

The above results distinguish, by a few simple properties, the locally
symmetric manifolds among all Riemannian manifolds.  We conjecture that
a stronger, more quantitative result holds, whereby there is a kind of
universal (depending only on $\pi_1$) constraint on the amount of
symmetry of any Riemannian manifold which is not an orbibundle with
locally symmetric fiber.

\begin{conjecture}
\label{conjecture:uniform}
The hypothesis ``$\Isom(\widetilde{M})$ is not discrete'' in
Theorem \ref{theorem:locsym1}, Corollary
\ref{corollary:locsym2}, and Corollary
\ref{corollary:locsym3} can be replaced by:
$[\Isom(\widetilde{M}):\pi_1(M)]>C$, where $C$ depends only on
$\pi_1(M)$.
\end{conjecture}

We do not know how to prove Conjecture \ref{conjecture:uniform}.
However, we can prove it in the special case of a fixed manifold
admitting a locally symmetric metric.

\begin{theorem}
\label{theorem:magic:number}
Let $(M,g_0)$ be a closed, irreducible, nonpositively curved 
locally symmetric $n$-manifold, $n>1$.  Then there exists a
constant $C$, depending only on $\pi_1(M)$, such that for any Riemannian
metric $h$ on $M$:
$$[\Isom(\widetilde{M}):\pi_1(M)]>C \mbox{\ \ if and only if\ \ \ }h\sim g_0$$
where $\sim$ denotes ``up to homothety of direct factors''.
\end{theorem}

\bigskip
\noindent
{\bf Manifolds with both closed and finite volume quotients. }
We can also apply our methods to answer the following
fundamental question in Riemannian geometry: which contractible 
Riemannian manifolds $X$ cover both a
closed manifold and a (noncompact, complete) finite volume
manifold?  

This question has been answered for many (but not all) contractible 
homogeneous spaces $X$.  Recall that a contractible 
{\em (Riemannian) homogeneous space} $X$ is the quotient
of a connected Lie group $H$ by a maximal compact subgroup, endowed with
a left-invariant metric.  Mostow proved that solvable $H$ admit only 
cocompact lattices, while 
Borel proved that noncompact, 
semisimple $H$ have both cocompact and noncocompact
lattices (see \cite{Ra}, Thms. 3.1, 14.1).  The case of arbitrary
homogeneous spaces is more subtle, and as far as we can tell remains
open.   

The following theorem extends these results to all contractible
manifolds $X$.  It basically states that if $X$ covers both a compact and
a noncompact, finite volume manifold, 
then the reason is that $X$ is ``essentially'' a product, with one
factor a homogeneous space which itself covers both types of manifolds.  

To state this precisely, we define 
a {\em warped Riemannian product} to be a smooth 
manifold $X=Y\times Z$ where $Z$ is a (locally) homogeneous space, 
$f:Y\to {\cal H}(Z)$ is a smooth function with target the space ${\cal
H}(Z)$ of all (locally) homogeneous metrics on $Z$, 
and the metric on $X$ is given by 
$$g_X(y,z)=g_Y\oplus f(y)g_Z$$

We can now state the following.

\begin{theorem}
\label{theorem:covering}
Let $X$ be a contractible Riemannian manifold.  
Suppose that $X$ Riemannian covers both a closed manifold and
a noncompact, finite volume, complete manifold.  Then $X$ is isometric 
to a warped product $Y\times X_0$, where $Y$ is a contractible
manifold (possibly a point) and $X_0$ is a 
homogeneous space which admits both cocompact and noncocompact lattices.
In particular, if $X$ is not a Riemannian warped product then
it is homogeneous.
\end{theorem}

Note that the factor $Y$ is necessary, as one can see by taking the
product of a homogeneous space with the universal cover of any compact
manifold.    We begin the deduction of Theorem \ref{theorem:covering} 
from the other results in this paper by noting that its hypotheses 
imply that $\Isom(Z)$ is nondiscrete, so that our general result can be 
applied.

\bigskip
\noindent
{\bf Irreducible lattices in products. }Let $X=Y\times Z$ be a
Riemannian product.  Except in obvious cases, $\Isom(Y)\times
\Isom(Z)\hookrightarrow \Isom(X)$ is a finite index inclusion.  
Recall that a lattice $\Gamma$ in $\Isom(X)$ is {\em irreducible} if it
is not virtually a product.  Understanding which Lie groups admit
irreducible lattices is a classical problem; see, e.g., \cite{Ma}, \S
IX.7.  Eberlein determined in \cite{Eb1,Eb2} the
nonpositively curved $X$ which admit irreducible lattices; they are
essentially the symmetric spaces.  The following extends this result to
all contractible manifolds; it also provides 
another proof of Eberlein's result.

\begin{theorem}
\label{theorem:irreducible}
Let $X$ be a nontrivial Riemannian product, and suppose that $\Isom(X)$
admits an irreducible, cocompact lattice.  Then $X$ is isometric to a warped
Riemannian product $X=Y\times X_0$, where $Y$ is a contractible
manifold (possibly a point), $X_0$ is a positive dimensional 
homogeneous space, and $X_0$ admits an irreducible, cocompact lattice.
\end{theorem}

As with Theorem \ref{theorem:covering}, Theorem
\ref{theorem:irreducible} is deduced 
from the other results in this paper by noting that its hypotheses 
imply that $\Isom(Z)$ is nondiscrete; see \S\ref{section:irreducible}.

\bigskip
\noindent
{\bf Compact complex manifolds. }Our results on isometries 
also have implications for complex manifolds.  Kazhdan conjectured that any
irreducible bounded domain $\Omega$ which admits both a compact quotient
$M$ and a one-parameter group of holomorphic automorphisms must be
biholomorphic to a bounded {\em symmetric} domain.  Frankel \cite{Fr1} 
first proved this for convex domains $\Omega$, and subsequent work by
Nadel \cite{Na} and Frankel \cite{Fr2}, which we now recall, proved it
in general.  

The Bergman volume form on a bounded domain produces a metric on the
canonical bundle so that the first Chern class satisfies $c_1(M)<0$;
equivalently, the canonical line bundle is ample.  
Hence Kazhdan's conjecture is implied by (and, indeed, inspired) 
the following.

\begin{theorem}[Nadel, Frankel]
\label{theorem:frankel}
Let $M$ be a compact, aspherical 
complex manifold with $c_1(M)<0$.  
Then there is a holomorphic splitting $M'=M_1\times M_2$ 
of a finite cover $M'$ of $M$, 
where $M_1$ is locally symmetric and $M_2$ is {\em locally
rigid} (i.e. the biholomorphic automorphism group of 
the universal cover $\widetilde{M_2}$ is discrete.)
\end{theorem}

Theorem \ref{theorem:frankel} was first 
proved in (complex) dimension two by Nadel \cite{Na} and in all
dimensions by Frankel \cite{Fr2}.  They do not
require the asphericity of $M$, although this is of course the case for  
quotients of bounded domains. Complex geometry is an essential 
ingredient in their work.

In \S\ref{section:complex} we give a different proof of Theorem
\ref{theorem:frankel}, using a key proposition from the earlier paper of
Nadel \cite{Na} together with (the proof of) our Theorem
\ref{theorem:main} below.  In complex dimension two, we give a proof 
independent of both \cite{Na} and \cite{Fr2}.  We do not see, however,
how to use our methods without the asphericity assumption.  

As with \cite{Na} and \cite{Fr2}, our starting point is a theorem of
Aubin-Yau, which gives that the biholomorphism group
$\Aut(\widetilde{M})$ acts isometrically on a Kahler-Einstein metric
lifted from $M$.  Our proof shows that, at least in complex dimension
two, this is the only ingredient from complex geometry needed to prove
Kazhdan's conjecture.

\bigskip
\noindent
{\bf Remark. }Nadel pointed out explicitly in Proposition 0.1 of
\cite{Na} that his methods would extend to prove Theorem \ref{theorem:frankel} 
in all dimensions if one could prove that 
each isotropy subgroup of $\Aut(\widetilde{M})^\circ$ were a
maximal compact subgroup.  The solution to this problem in the 
aspherical case is given in Claim IV of Section \ref{section:proof} below;
it also applies outside of the holomorphic context as well.  

\bigskip
\noindent
{\bf Some additional applications. }
A number of the results from this paper generalize from closed, aspherical 
Riemannian manifolds to all closed Riemannian manifolds.  
In \S\ref{section:nonaspherical} we provide an illustrative example,  
Theorem \ref{theorem:nonaspherical}, which seems to be the
first geometric rigidity theorem for non-aspherical manifolds with
infinite fundamental group.

In \S\ref{section:Hopf} below we give an application 
of our methods to the Hopf Conjecture about Euler characetristics 
of aspherical manifolds.  

Finally, we would like to mention the work of K.\ Melnick in
\cite{Me}, where some of the results here are extended from the
Riemannian to the pseudo-Riemannian (especially the Lorentz) case.  
Melnick combines the ideas here with Gromov's theory of rigid geometric
structures, as well as methods from Lorentz dynamics.

\bigskip
\noindent
{\bf Acknowledgements. }A first version of the main results of this
paper were proved in the Fall of 2002.  We would like to thank the
audiences of the many talks we have given since that time on the work
presented here; they provided numerous useful comments.  
We are particularly grateful to the students in
``Geometric Literacy'' at the University of Chicago, especially to Karin
Melnick for her corrections on an earlier version of this paper.  We 
would like to thank Ralf Spatzier who, after hearing a talk on some
of our initial results (later presented in
\cite{FW}), pointed out a connection with Eberlein's work; 
this in turn lead us to the idea that a much more general
result might hold.  Finally, we would like to thank the excellent referees, 
whose extensive comments and suggestions greatly improved the paper.

\section{Finding the orbibundle (proof of Theorem \ref{theorem:main})}
\label{section:proof}

Our goal in this section is to prove Theorem \ref{theorem:main}.  
The starting point is the following well-known classical theorem.

\begin{theorem}[Myers-Steenrod, \cite{MS}]
\label{theorem:myers-steenrod}
Let $M$ be a Riemannian manifold.  Then $\Isom(M)$ is a Lie group, and
acts properly on $M$.  If
$M$ is compact then $\Isom(M)$ is compact.
\end{theorem}

Note that the Lie group $\Isom(M)$ in Theorem
\ref{theorem:myers-steenrod} may have infinitely many components; for
example, let $M$ be the universal cover of a bumpy metric on
the torus.

Throughout this paper we will use the following notation:

\bigskip
\noindent
$M$= a closed, aspherical Riemannian manifold

\smallskip
\noindent
$\Gamma=\pi_1(M)$

\smallskip
\noindent
$X=\widetilde{M}=$ the universal cover of $M$ 

\smallskip
\noindent
$I=\Isom(X)=$ the group of isometries of $X$

\smallskip
\noindent
$I_0=$ the connected component of $I$ containing the identity

\smallskip
\noindent
$\Gamma_0=\Gamma\cap I_0$

\bigskip
Here $X$ is endowed with 
the unique Riemannian metric for which the covering map $X\rightarrow M$
is a Riemannian covering.    
Hence $\Gamma$ acts on $X$ isometrically 
by deck transformations, giving a natural inclusion $\Gamma\rightarrow
I$, where $I=\Isom(X)$ is the isometry group of $X$.  

By Theorem \ref{theorem:myers-steenrod}, $I$ is a Lie group, possibly
with infinitely many components.  Let $I_0$ denote the connected
component of the identity of $I$; note that $I_0$ is normal in $I$.  
If $I$ is discrete, then we are done,
so suppose that $I$ is not discrete.  Theorem
\ref{theorem:myers-steenrod} then gives that the dimension of $I$ is
positive, and so $I_0$ is a connected, 
positive-dimensional Lie group. 

We have the following exact sequences:

\begin{equation}
\label{eq:exact1}
1\to I_0\to I\to I/I_0\to 1
\end{equation}

and 
\begin{equation}
\label{eq:exact2}
1\to \Gamma_0\to \Gamma\to \Gamma/\Gamma_0\to 1
\end{equation}

We now proceed in a series of steps.  Our first step is to 
construct what will end up as the locally homogeneous fibers of the
orbibundle (\ref{eq:general:orbi}).  

\bigskip
\noindent
{\bf Claim I: }The quotient $I_0/\Gamma_0$ is compact.

\proof
Let $\frames(X)$ denote the frame bundle over $X$.  The
isometry group $I$ acts freely on $\frames(X)$. The $I_0$ orbits
in $\frames(X)$ give a smooth foliation of $\frames(X)$ whose
leaves are diffeomorphic to $I_0$.  This foliation descends via
the natural projection $\frames(X)\to \frames(M)$ to give a smooth
foliation ${\cal F}$ on $\frames(M)$, each of whose leaves is
diffeomorphic to $I_0/\Gamma_0$.  Thus we must prove that each of 
these leaves is compact.

The quotient of $\frames(X)$ by the smallest subgroup of $I$ 
containing both $\Gamma$ and $I_0$ is homeomorphic to the space of leaves of
${\cal F}$.  We claim that this quotient is a finite cover 
of $\frames(X)/I$.  To prove this, it is clearly enough to show that the
natural injection $\Gamma/\Gamma_0\to I/I_0$ has finite index
image.  

To this end, we first recall the following basic principle of 
Milnor-Svarc (see, e.g., \cite{H}).  
Let $G$ be a compactly generated topological group, 
generated by a compact subspace $S\subset G$.  Endow $G$ with the 
{\em word metric}, i.e. let $d_G(g,h)$ be defined to be the minimal
number of elements of $S$ needed to represent $gh^{-1}$; this is a 
left-invariant metric on $G$.  Now suppose that $G$ acts properly and
cocompactly by isometries on a proper, geodesic metric space $X$.
Then $G$ is {\em quasi-isometric} to $X$, i.e. for any fixed basepoint 
$x_0\in X$, the orbit map $G\to X$
sending $g$ to $g\cdot x_0$ satisfies the following two conditions: 
\begin{itemize}
\item (Coarse Lipschitz): For some $K,C>0$, 
$$\frac{\displaystyle 1}{\displaystyle K}d_G(g,h)-C\leq d_X(g\cdot
x_0,h\cdot x_0)\leq Kd_G(g,h)+C$$
\item($C$-density) $\nbhd_C(G\cdot x)=X$
\end{itemize}

While the standard proofs of this fact (see, 
e.g., \cite{H}) usually assume that $S$ is
finite, they apply verbatim to the more general case of $S$ compact.  

Applying this principle, the cocompactness of the actions of both 
$\Gamma$ and of $I$ on $X$ give that 
the inclusion $\Gamma\to I$ is a quasi-isometry.  The quotient map 
$I\to I/I_0$ is clearly distance nonincreasing, and so the image
$\Gamma/\Gamma_0$ of $\Gamma$ under this quotient map is $C$-dense 
in $I/I_0$.  As both groups are discrete, this clearly implies that the 
inclusion $\Gamma/\Gamma_0\to I/I_0$ is of finite index.
Thus the claim is proved.

Now note that $\frames(X)/I$ is clearly compact, and is a manifold since
$I$ is acting freely and properly.  Hence the leaf-space of ${\cal F}$
is also a compact manifold.  Since each leaf of ${\cal F}$ is the
inverse image of a point under the map from $\frames(M)$ to the leaf
space, we have that each leaf of ${\cal F}$ is compact.
\endproof

It will be useful to know that $I_0$ cannot have compact factors.

\bigskip
\noindent
{\bf Claim II: }$I_0$ has no nontrivial compact factor.

\proof 
In proving this claim, we will use degree theory for noncompact
manifolds, phrased in terms of locally finite homology $H_\ast^{lf}$
(see, e.g., \cite{Iv} for a discussion).  
{\em Locally finite homology} is the theory 
of cycles which pair with cohomology with compact support.  
Perhaps the quickest description of $H_\ast^{lf}(X)$ is as
the usual reduced homology $\widetilde{H}_\ast(\widehat{X})$ of the one-point
compactification $\widehat{X}$ of $X$.  Alternatively, it can be
described (for locally finite simplicial complexes) as the
homology of the chain complex of infinite formal combinations of
simplices for which only finitely many simplices with nonzero
coefficients intersect any given compact region.

With this definition, it is easliy verified (see \cite{Iv})
that the usual degree theory holds for continuous quasi-isometries
between (possibly noncompact) manifolds $X$, with the fundamental class 
of the $n$-manifold $X$ now being an element of $H_n^{lf}(X,\Z)$.  
As one example, the universal cover $X$ of a closed, aspherical
$n$-manifold $M$ has a nonzero fundamental class lying in
$H_n^{lf}(X,\Z)$.  With this degree theory in place, we can now 
begin the proof of Claim II.
 
Now suppose that $I_0$ has a nontrivial compact factor $K$.  Since $I_0$
is connected and $\dim(I_0)>0$ by assumption, we have that $K$ is
connected and $\dim(K)>0$.
  
Since $M$ is closed, and so there is a compact 
fundamental domain for the $\Gamma$-action on $X$, 
we easily see that there exists a
constant $C$ so that each $K$-orbit has diameter at most $C$.  But
then $X/K$ is quasi-isometric to $X$.  Now the standard
``connect the dots'' trick (see, e.g., p.527 of \cite{BW}, or 
Appendix A of \cite{BF} for exact
details) states that such quasi-isometries are a bounded distance (in the
$\sup$ norm) from a {\em continuous} quasi-isometry (i.e. Lipschitz
map).  Hence 
there are continuous maps $X\to X/K$ and $X/K\to X$
inducing the given quasi-isometry.    
Since $\dim(K)>0$ we have that 
$\dim(X/K)<\dim(X)=n$.  This implies that the fundamental
class of $X$ in $H^{lf}_n(X,\R)$, where $n=\dim(X)$,
must vanish, contradicting the fact that $X$ is the universal cover of a
closed, aspherical $n$-manifold.
\endproof

The next step in our proof of Theorem \ref{theorem:main}
is to determine information which will help us
construct the orbifold base space $B$ 
of the orbibundle (\ref{eq:general:orbi}).

\bigskip
\noindent
{\bf Claim III: }$X/I_0$ is contractible.

\proof  The Conner Conjecture,
proved by Oliver \cite{Ol}, gives that the quotient of
a contractible manifold by a connected, compact, smooth transformation 
group is contractible.  Our claim that $X/I_0$ is contractible 
follows directly from the following simple extension of Oliver's theorem.  

\begin{proposition}
\label{theorem:conner}
Let $G$ be a connected Lie group acting properly by diffeomorphisms on a
contractible manifold $X$.  Then the underlying topological space of the
orbit space $X/G$ 
is contractible.
\end{proposition}

Proposition \ref{theorem:conner} is a consequence of Oliver's Theorem
and the following.

\begin{proposition}
Let $G$ be a connected Lie group acting properly by diffeomorphisms on
an aspherical manifold $X$.  Denote by $K$ the maximal 
compact subgroup of $G$. Then there exists an aspherical 
manifold $Y$ such that $X$ is diffeomorphic to $Y\times G/K$, the
manifold $Y$ has a $K$-action, and the original action is given by the 
product action.  In particular, $X/G$ is diffeomorphic to $Y/K$. 
\end{proposition}

\proof
Let $\underline{EG}$ be the classifying space for proper
CW $G$-complexes, so that $\underline{EG}/G$ is the classifying space
for proper $G$-bundles (see, e.g., the appendix of \cite{BCH}).   
Now $G/K$ is an $\underline{EG}$ space.  Hence 
there is a proper $G$-map
$\psi:X\to \underline{EG}$.  But $\underline{EG}$ has 
only one $G$-orbit, so $\psi$ is
surjective.  Now let $Y=\psi^{-1}([K])$, where $[K]$ denotes the identity
coset of $K$. Hence $X$ is diffeomorphic to $G\times_KY$, and we
are done.

\endproof

We are now ready to construct, on the level of universal covers, the 
orbibundle (\ref{eq:general:orbi}), and in particular to prove that the
base space $B$ is a Riemannian orbifold.  The crucial point is to understand
stabilizers of the $I_0$ action on $X$.  
For $x\in X$, denote the stabilizer of $x$ under the $I_0$ action by 
$I_x:=\{g\in I_0: gx=x\}$.  Let $K_0$ denote the maximal compact
subgroup of $I_0$; it is unique up to conjugacy.  

\bigskip
\noindent
{\bf Claim IV: }$I_x=K_0$ for each $x\in X$.  Hence the following hold:
\begin{enumerate}
\item $X/I_0$ is a manifold. 
\item Each $I_0$-orbit in $X$ is isometric to the 
contractible, homogeneous manifold $I_0/K_0$, endowed with some
left-invariant Riemannian metric.
\item The natural quotient map gives a Riemannian fibration
\begin{equation}
\label{eq:fibration1}
I_0/K_0\to X\to X/I_0
\end{equation}
\end{enumerate}

\proof  
Clearly $I_x\subseteq K_0$.  
Iwasawa proved (\cite{Iw}, Theorem 6) 
that any maximal compact subgroup of a connected Lie group
is connected.  Hence it is enough to prove that
$\dim(I_x)=\dim(K_0)$.  

To this end we consider rational cohomological dimension $\cdq$.  By 
Claim III we have $X/I_0$ is contractible.  
Since $\Gamma/\Gamma_0$ acts properly on
$X/I_0$, we then have
\begin{equation}
\label{eq:cd1}
\cdq(\Gamma/\Gamma_0)\leq \dim(X/I_0)
\end{equation}
Since $K_0$ is maximal, we know $I_0/K_0$ is contractible.  By Claim I,
we have that $\Gamma_0$ is a uniform lattice in $I_0$, and so
\begin{equation}
\label{eq:cd2}
\cdq(\Gamma_0)=\dim(I_0/K_0)
\end{equation} 

Since $X$ is contractible and $M=X/\Gamma$ is a closed manifold, 
by general facts about cohomological dimension (see \cite{Bro}, Chapter
VIII (2.4)), we have
$$\dim(X)=\cdq(\Gamma)\leq \cdq(\Gamma_0)+\cdq(\Gamma/\Gamma_0)$$
which combined with (\ref{eq:cd1}) and (\ref{eq:cd2}) gives

\begin{equation}
\label{eq:cd3}
\dim(X)\leq \dim(X/I_0)+\dim(I_0/K_0)
\end{equation}

But for each $x\in X$, we have
\begin{equation}
\label{eq:cd4}
\dim(X)\geq \dim(X/I_0)+\dim(I_0/I_x)
\end{equation}
which combined with (\ref{eq:cd3}) gives 
$$\dim(I_0/I_x)\leq \dim(I_0/K_0)$$
and so $\dim(I_x)\geq \dim(K_0)$, as desired.   Thus $I_x=K_0$.  

It follows that each orbit $I_0\cdot x$ is diffeomorphic to a common
Euclidean space $I_0/K_0$, so by the Slice Theorem 
(see, e.g., \cite{Br}, Chap. IV, \S 3,4,5) 
it follows that $X/I_0$ is a manifold.  We note that while 
the Slice Theorem is
usually stated for actions of compact groups, the proof extends 
immediately to the case of proper actions of noncompact groups; one
simply produces an invariant Riemannian
metric by translating a compactly supported pseudometric, and this gives
the required structure via exponentiation.
\endproof

\bigskip
\noindent
{\bf Finishing the proof. }
The action of $\Gamma$ on $X$ induces actions of $\Gamma_0$ on
$I_0/K_0$, and of $\Gamma/\Gamma_0$ on $X/I_0$, compatible with the
Riemannian fibration (\ref{eq:fibration1}).  By Myers-Steenrod, 
$\Gamma/\Gamma_0$ acts properly discontinuously on
$X/I_0$; we denote the quotient space of this action by $B$.  
We thus have a Riemannian orbibundle
(as defined in the introduction): 

\begin{equation}
\label{eq:fibration2}
F\to M\to B
\end{equation}
where $F$ denotes the closed, locally homogeneous Riemannian manifold 
$\Gamma_0\backslash I_0/K_0$, endowed with the quotient metric of 
a left $I_0$-invariant metric on $I_0/K_0$.  This completes the proof of
Theorem \ref{theorem:main}

\section{The case when \boldmath$I_0$ is semisimple}

The main goal of this section is to prove Proposition
\ref{proof:semisimp:consequence} below, 
which shows that when $I_0$ is semisimple with finite center, a
much stronger conclusion holds in Theorem \ref{theorem:main}.  

\begin{proposition}
\label{proof:semisimp:consequence}
Suppose that $I_0$ is semisimple with finite center.  Then $M$ has a
finite cover which is a Riemannian warped product $N\times B$, where $N$
is nonempty, locally symmetric with nonpositive curvature, and has 
no local torus factors.  In particular, $\pi_1(B)\lhd\pi_1(M)$, and 
any nontrivial, normal abelian subgroup of $\pi_1(M)$ lies in $\pi_1(B)$. 
\end{proposition}

\medskip
\noindent
{\bf Remark on semisimplicity. }We would like to emphasize that by
calling a connected Lie group $G$ ``semisimple'' we mean only that the
Lie algebra of $G$ is semisimple.  Thus the center $Z(G)$ may be
infinite.  Such examples do exist (for example the universal cover of
$\U(n,1)$) , and must be taken into account.  We
also point out that $G$ may in general have compact factors.  For
the connected component $I_0$ of the isometry group of the universal
cover of a closed, aspherical 
Riemannian manifold, however, we have already proven 
in Claim II of \S\ref{section:proof} above that $I_0$ has no nontrivial
compact factor.  Even so, the semisimple part $(I_0)^{\sem}$ may have 
nontrivial compact factors coming from $Z(I_0)$.  
\medskip

After proving Proposition \ref{proof:semisimp:consequence}, 
we show in 
\S\ref{section:nonormal:proof} that the hypothesis that $I_0$ is
semisimple with finite center is more common than one might guess.
Indeed, in Proposition \ref{prop:no:normal:abelian} we prove that $I_0$
is always semisimple with finite center unless $\Gamma=\pi_1(M)$
contains an infinite, normal abelian subgroup.   

\subsection{The proof of Proposition \ref{proof:semisimp:consequence}}
\label{section:when:semisimp}

The structure of the proof of Proposition
\ref{proof:semisimp:consequence} is to first prove it at the level of
fundamental groups, mostly using Lie theory.  The theory of harmonic 
maps, as well as the existence of arithmetic lattices, is then used to 
build so many isometries of the universal cover of $M$ that it is forced
to fiber in the claimed way.

\bigskip
\noindent
{\bf Triviality of the extension. }
Our first goal will be to prove that, after replacing $\Gamma$ by a
finite index subgroup if necessary, the exact sequence 

\begin{equation}
\label{eq:exact5}
1\to \Gamma_0\to \Gamma\to \Gamma/\Gamma_0\to 1
\end{equation}
splits as a direct product.  As with every extension, (\ref{eq:exact5}) is
determined by two pieces of data: 
\begin{enumerate}
\item A representation $\rho:\Gamma/\Gamma_0\to \Out(\Gamma_0)$, and 
\item A cohomology class in $H^2(\Gamma/\Gamma_0,Z(\Gamma_0)_\rho)$,
where $Z(\Gamma_0)_\rho$ is a $\Gamma/\Gamma_0$-module via $\rho$.
\end{enumerate}

We analyze these pieces in turns.  Let $<I_0,\Gamma>$ be the smallest
subgroup of $I$ containing $I_0$ and $\Gamma$.  Consider the exact
sequence

\begin{equation}
\label{eq:ss4}
1\to I_0\to <I_0,\Gamma>\to \Gamma/\Gamma_0\to 1
\end{equation}

and let
$$\rho_1:\Gamma/\Gamma_0\to \Out(I_0)$$
denote the induced action; this is just the action induced by the
conjugation action of $\Gamma$ on $I$.  
Since $I_0$ is semisimple, we know (see, e.g. \cite{He}, Theorem IX.5.4)
that $\Out(I_0)$ is finite.
Hence, after passing to a finite index subgroup of $\Gamma$ if
necessary, we may assume that $\rho_1$ is trivial.  In other words, the
$\Gamma$-action on $I_0$ is by inner automorphisms, giving a
representation 
$$\rho_2:\Gamma/\Gamma_0\to I_0/Z(I_0)$$

Now, the conjugation action of $\Gamma$ on $I_0$ preserves $\Gamma_0$,
and so the image of $\rho_2$ lies in the normalizer $N_H(\Gamma_0)$ of
$\Gamma_0$ in $H:=I_0/Z(I_0)$. Note that $\Gamma_0\cap Z(I_0)$ is finite
and hence trivial, as is $Z(\Gamma_0)$, 
since $\Gamma_0$ is torsion free, and so $\Gamma_0$
can be viewed as a subgroup of $H$.  Since $H$ is semisimple and
$\Gamma_0$ is a cocompact lattice in $H$ (by Claim I in the
proof of Theorem \ref{theorem:main}), it follows that 
$N_H(\Gamma_0)/\Gamma_0$ is finite \footnote{This follows for example from
Bochner's classical result that the closed manifold 
$M=\Gamma\backslash H/K$ has finite isometry group since it has negative
Ricci curvature, and $\Isom(M)=N_H(\Gamma_0)$.  For another proof, see 
\cite{Ma}, II.6.3.}.

Hence, by replacing $\Gamma$ with a finite index subgroup if necessary, 
we may assume $\rho_2$ has trivial image.  We thus have that the
conjugation action of $\Gamma$ on $\Gamma_0$ is by inner automorphisms 
of $\Gamma_0$.  Since $Z(\Gamma_0)$ is trivial, 
the representation $\rho:\Gamma/\Gamma_0\to \Out(\Gamma_0)$ 
is trivial \footnote{Note that there are cases when $\Out(\Gamma_0)$ is
nontrivial; for example when $\Gamma_0$ is a surface group then
$\Out(\Gamma_0)$ is the mapping class group of that surface.}.  We also
know that 

$$H^2(\Gamma/\Gamma_0,Z(\Gamma_0)_{\rho})=0$$
since $Z(\Gamma_0)=0$.  It follows that, up to finite index, the exact sequence
(\ref{eq:exact5}) splits, and in fact that

\begin{equation}
\label{eq:groups:split}
\Gamma\approx \Gamma_0\times \Gamma/\Gamma_0
\end{equation}

Recall (\ref{eq:fibration2}), where we found a Riemannian orbibundle 
$$F\to M\to B$$
Our goal now is to use (\ref{eq:groups:split}) to find a section of this
fibration, and to use this to prove that $M$ is a Riemannian warped
product.  In order to do this we will use the following tool.

\bigskip
\noindent
{\bf Harmonic maps. }
We recall that a map $f:N\to M$ between Riemannian manifolds is {\em
harmonic} if it minimizes the energy functional
$$E(f)=\int_N||Df_x||^2d{\rm vol}_N$$

The key properties of harmonic maps between closed Riemannian 
manifolds which we will need are the following (see, e.g. \cite{SY}): 

\begin{itemize}
\item (Eels-Sampson) 
When the target manifold has nonpositive sectional curvatures, a
harmonic map exists in each homotopy class.

\item (Hartman, Schoen-Yau) 
If a harmonic map $f:M\to N$ induces a surjection on $\pi_1$,
and if $\pi_1(N)$ is centerless, then $f$ is unique in its homotopy
class. This follows directly from Theorem 2 of \cite{SY}.

\item (easy) The precomposition and postcomposition of 
a harmonic map with an isometry gives a harmonic map.

\end{itemize}

\bigskip
\noindent
{\bf Showing that $X$ is a warped product. } The isomorphism in
(\ref{eq:groups:split}) gives via projection to a direct factor a
natural surjective homomorphism $\pi:\Gamma\to\Gamma_0$.  Recall that 
$0=Z(\Gamma_0)\supseteq Z(I_0)\cap \Gamma_0$, and so the injection
$\Gamma_0\to I_0$ gives an injection $\Gamma_0\to I_0/Z(I_0)$.  Our 
first goal is to extend the projection $\pi$ to a projection 
$\widehat{\pi}:<I_0,\Gamma>\to I_0/Z(I_0)$.

To this end, note
that $Z(I_0)$ is characteristic in $I_0$, and so $Z(I_0)\lhd I$; in
particular, $Z(I_0)\lhd <I_0,\Gamma>$.  Taking the 
quotient of the exact sequence (\ref{eq:ss4}) by the finite normal 
subgroup $Z(I_0)$ gives an exact sequence
\begin{equation}
\label{eq:ss4b}
1\to I_0/Z(I_0)\to <I_0,\Gamma>/Z(I_0)\to \Gamma/\Gamma_0\to 1
\end{equation}

We claim that the kernel of (\ref{eq:ss4b}) is centerless.  Indeed, if
$G$ is any connected semisimple Lie group, then its center $Z(G)$ is
clearly closed, hence discrete since $G$ is semsimple. But for any connected
Lie group $G$ with $Z(G)$ discrete, the center of $G/Z(G)$ is trivial 
(see, e.g., Exercise 7.11(b) of \cite{FH})\label{s4}.  The reason this fact is
true can be seen from the fact that the discreteness of $Z(G)$ implies 
that the quotient map $G\to G/Z(G)$ is a covering map of Lie groups, 
and so both $G$ and $G/Z(G)$ have isomorphic Lie algebras and isomorphic
universal covers..  

Since the kernel of (\ref{eq:ss4b}) is centerless, the exact argument as
above gives that (\ref{eq:ss4b}) splits, so that 
\begin{equation}
\label{eq:ss4c}
<I_0,\Gamma>/Z(I_0)\approx I_0/Z(I_0)\times \Gamma/\Gamma_0 
\end{equation}

This isomorphism, composed with the natural projections, then gives us a
surjective homomorphism $$\widehat{\pi}:<I_0,\Gamma>\to I_0/Z(I_0)$$

Let $K_0$ denote a maximal compact subgroup of the semisimple Lie group 
$I_0$.  We then 
have that $I_0$ acts isometrically on the contractible, nonpositively
curved symmetric space of noncompact type 
$X_0:=I_0/K_0$.  Since $Z(I_0)$ is finite, it
lies in $K_0$, and so the $I_0$ action on $X_0$ factors through a  
faithful action of $I_0/Z(I_0)$.  As $X_0$ is contractible, the homomorphism 
$\pi$ is induced by some continuous map $h:X/\Gamma\to X_0/\Gamma_0$.  
Thus $f$ is homotopic to a harmonic map $h$.  By the theorem of Hartman
and Schoen-Yau stated above, $f$ is the unique harmonic map in its homotopy
class.

\begin{claim}
\label{claim:equiv}
The lifted map $\widetilde{f}:X\to X_0$ is equivariant with respect to
the representation $\widehat{\pi}:<I_0,\Gamma>\to I_0/Z(I_0)$.
\end{claim}

To prove this claim, first note that $\widetilde{f}$ is equivariant with
respect to the representation $\pi$, by construction; we want
to promote this to $\widehat{\pi}$-equivariance.  One strange aspect of
this is that we use an auxilliary arithmetic group, which seems to have 
nothing to do with the situation.

To begin, consider any cocompact lattice $\Delta$ in $I_0/Z(I_0)$.  
By (\ref{eq:ss4c}), $\Delta\times \Gamma/\Gamma_0$ is a cocompact 
lattice in $<I_0,\Gamma>/Z(I_0)$, so it pulls back 
under the natural quotient to a cocompact 
lattice, which we will also call $\Delta$, in $<I_0,\Gamma>$ (recall
that $Z(I_0)$ is finite).  
  
Then, up to translation by elements of
$\Delta$, there is a unique harmonic map 
$\phi_\Delta:X\to X_0$ equivariant with respect to the restriction 
$\widehat{\pi}|_{\Delta\times (\Gamma/\Gamma_0)}$.  Suppose $\Delta'$ is any other lattice in $I_0$
which is commensurable with $\Delta$.  Since both $\phi_\Delta$ and 
$\phi_{\Delta'}$ are harmonic and equivariant with respect to the 
representation $\widehat{\pi}$ restricted to 
$(\Delta\cap \Delta')\times (\Gamma/\Gamma_0)$, and
since $\Delta\cap \Delta'$ has finite index in both $\Delta$ and in 
$\Delta'$, we have by uniqueness of harmonic maps that
$\phi_\Delta=\phi_{\Delta'}$.  We remark that this ``uniqueness implies
equivariance'' principle is also a key trick in \cite{FW}.  

Since $I_0$ is semisimple with finite center, the quotient $I_0/Z(I_0)$ is
semisimple and centerless (as proven just after equation ({eq:ss4b}) 
on page \pageref{s4} above).  By a theorem of Borel (\cite{Bo}, Theorem C), 
there exists a cocompact arithmetic 
lattice $\Delta_1$ in $I_0/Z(I_0)$. Since $I_0/Z(I_0)$ is centerless, 
it follows that the commensurator $\Comm_{I_0/Z(I_0)}(\Delta_1)$ is
dense in $I_0/Z(I_0)$; see, for example, 
Proposition 6.2.4 of \cite{Zi}, where this is clearly explained.  

Let $\Delta_0$ denote the pullback of $\Delta_1$ under the natural 
quotient map $I_0\to I_0/Z(I_0)$.  Since $\Delta_0$ contains $Z(I_0)$, 
and so $\Comm_{I_0}(\Delta_0)$ is the central extension of
$\Comm_{I_0/Z(I_0)}(\Delta_1)$ associated to $Z(I_0)$, it follows that
$\Comm_{I_0}(\Delta_0)$ is dense in $I_0$.  
At this point, a verbatim application of the proof of the
``arithmetic case'' of Theorem 1.4 in \cite{FW} completes the proof of 
Claim \ref{claim:equiv}; for completeness, we briefly recall this proof.
 
Let $U$ denote the set of $g\in <I_0,\Gamma>$ for which the
equation
\begin{equation}
\label{eq:equiv1}
\phi_{\Delta_0}g=g\phi_{\Delta_0}
\end{equation}
holds.  Now $U$ is closed, and the uniqueness of 
harmonic maps gives that $U$ is a subgroup of $I_0$.  Hence $U$ is a Lie
subgroup of $I_0$.  Applying
the above paragraphs with $\Delta=\Delta_0$ and with 
$\Delta'$ running through the collection $\cal L$ 
of lattices commensurable with $\Delta_0$ in $I_0$, gives that 
$U$ contains every lattice in $\cal L$.  Since 
$\Comm_{I_o}(\Delta_0)$ is dense in $I_0$, there are infinitely
many distinct members of $\cal L$ conjugate to $\Delta_0$, namely the
conjugates of $\Delta_0$ by elements of $\Comm_{I_o}(\Delta_0)$.  Hence 
$U$ is nondiscrete, hence
positive dimensional.  Under the adjoint representation, $\Delta_0$
preserves the Lie algebra of $U$.  But $\Delta_0$ is a lattice in $I_0$,
hence is Zariski dense by the Borel Density Theorem (see,
e.g. \cite{Ma}, Theorem II.2.5).  Thus $U=I_0$, finishing the proof 
of Claim \ref{claim:equiv}.

We now have a map 
$$X/(\Gamma \times (\Gamma/\Gamma_0))\to (X_0/\Gamma_0) \times
(X/I_0)/(\Gamma/\Gamma_0)$$
given by the product of $\widetilde{f}$ and the natural orbit map.  
This map harmonic when composed with projection to the first factor, and
is clearly a diffeomorphism, since we have just shown that the first
coordinate is equivariant with respect to $\widehat{\pi}$.

\subsection{Consequences of no normal abelian subgroups}
\label{section:nonormal:proof}

The assumption that $\Gamma=\pi_1(M)$ contains no nontrivial normal abelian
subgroup has strong consequences for our setup.  The main one is the
following.

\begin{proposition}
\label{prop:no:normal:abelian}
Suppose $\Gamma$ contains no infinite, normal abelian subgroup.  Then
$I_0$ is semisimple with finite center.
\end{proposition}

\proof
Note that since $\Gamma$ is torsion free, it follows that $\Gamma$ has
no normal abelian subgroups; in particular $Z(\Gamma)=1$.

For any connected Lie group $G$ there is an exact sequence 
\begin{equation}
\label{equation:gen1}
1\to G^{\sol}\to G\to G^{\sem} \to 1
\end{equation}
where $G^{\sol}$ denotes the {\em solvable radical} of $G$ (i.e. the maximal
connected, normal, solvable Lie subgroup of $G$), and where $G^{\sem}$ is
the connected semisimple Lie group $G/G^{\sol}$.  

Let $\Gamma^{\sol}$ denote the maximal normal solvable subgroup of
$\Gamma$; it is of course torsion free since $\Gamma$ is torsion free.
We claim that $\Gamma^{\sol}$ is trivial.  Suppose not.  Being a
nontrivial torsion-free solvable group, $\Gamma^{\sol}$ would then have an
infinite, characteristic, torsion-free abelian subgroup $H$, namely the
last nontrivial term in its derived series.  Since $H$ is characteristic
in the normal subgroup $\Gamma^{\sol}$ of $\Gamma$, it would follow that $H$ is
normal in $\Gamma$.  Since $H$ is infinite abelian, this contradicts the
hypothesis on $\Gamma$.

We now quote a result of Prasad, namely Lemma 6 in
\cite{Pr}.  For a lattice $\Gamma$ in a Lie group $I_0$, 
Conclusion (2) of Prasad's Lemma gives, in the terminology of \cite{Pr}:
$$\rank(\Gamma^{\sol})=\chi(I_0^{\sol})+\rank(Z(I_0^{\sem}))$$ Here
$\chi(I_0^{\sol})$ denotes the dimension of $I_0^{\sol}$ minus that of
its maximal compact subgroup, $\rank(Z(I_0^{\sem}))$ denotes the the rank of
the center of $I_0^{\sem}$, and {\em rank} denotes the sum of the ranks
of the abelian quotients in the derived series.  Since in our case we
have proven that $\Gamma^{\sol}=0$, it follows both that
$\chi(I_0^{\sol})=0$, i.e. that $I_0^{\sol}$ is compact, and that the
rank of $Z(I_0)$ is $0$, so that $Z(I_0)$ is finite.

Since $I_0^{\sol}$ is both solvable and compact, it is a torus $T$.  
Since the automorphism group of $T$ is discrete (namely it is
$\GL(\dim(T),\Z)$), the natural conjugation action of the 
connected group $I_0^{\sem}$ on $T$ given by (\ref{equation:gen1}) 
must be trivial, so that $T$ is a direct factor of $I_0$.  But we have
already proven (Claim II of \S\ref{section:proof}) 
that $I_0$ has no nontrivial compact factors, a contradiction unless 
$T$ is trivial.  
Thus $I_0^{\sol}=T$ is trivial; that is, $I_0$ is semisimple.
\endproof

\medskip
\noindent
{\bf Remark. }It is possible to weaken the hypothesis of Proposition 
\ref{prop:no:normal:abelian}, and hence of all of the results which rely
on it, to assuming only that $\Gamma$ contains no 
{\em finitely generated}, infinite normal abelian subgroups.  To do
this, we begin by recalling that Prasad's result used above also gives
that the group $\Gamma^{\sol}$ is a lattice in
some connected solvable subgroup $S$ of $I_0$.  
It follows from Proposition \ref{prop:mostow} below that $\Gamma^{\sol}$
is polycyclic.  But it is well-known and easy to see that 
any polycyclic group has the property that each of its
subgroups is finitely generated (see, e.g. \cite{Ra}, Prop. 3.8). 
Hence the subgroup $H$ constructed in the proof of Proposition 
\ref{prop:no:normal:abelian} would in fact be finitely generated.
\medskip

In the argument just given we needed the 
following proposition, proved by Mostow in the simply connected case.  

\begin{proposition}
\label{prop:mostow}
Every lattice $\Lambda$ in a connected solvable Lie group $S$ 
is polycyclic.
\end{proposition}

\proof  First note that $\pi_1(S)$ is
finitely-generated and abelian, and so the 
universal cover $\widetilde{S}$ is a central $\Z^d$ 
extension of $S$ for some $d\geq 0$.  The lattice 
$\Lambda$ pulls back to a lattice
$\widetilde{\Lambda}$ in $\widetilde{S}$, which is a central $\Z^d$ 
extension of $\Lambda$.  Mostow proved (see, e.g. \cite{Ra}, Prop. 3.7)
that any lattice $\widetilde{\Lambda}$ in a 
connected, simply-connected solvable Lie group $\widetilde{S}$ must be 
polycyclic. It follows easily that $\Lambda$ is polycyclic. 
\endproof

The use of Prasad's result simplifies the approach to Proposition 
\ref{prop:no:normal:abelian} given in an earlier version of this paper.
As part of that earlier approach, we proved the following proposition.
We include this result here since we believe it might prove useful in
the future, since the proof is direct, and since 
we were not able to find this result in the literature.  The argument
was kindly supplied to us by the referee.

\begin{proposition}
\label{prop:nonlinear1}
Let $G$ be a connected semsimple Lie group, and let $\Lambda$ be a 
lattice in $G$.  If $Z(G)$ is infinite then $Z(\Lambda)$ 
is infinite.
\end{proposition}

\proof
Let $T$ be the identity component of the closure of 
$Z(G)\Lambda$ in $G$.  First note that $T$ is abelian; indeed, 
the commutator subgroup $[T,T]$ of $T$ 
is contained in the closure of the subgroup
$$[Z(G)\Lambda, Z(G)\Lambda] = [\Lambda, \Lambda]\subset \Lambda$$ 
and hence $[T,T]$, being connected, is trivial.

Now let $C$ be the unique maximal compact, connected normal
subgroup of $G$. Then the Borel Density Theorem applied to the image
of $Z(G)\Lambda$ in $G/C$ gives that the image of $T$ in $G/C$ is a
connected, normal abelian subgroup.  Hence it must be trivial. 
Thus $T\subseteq C$, and so it is a torus normalized by $Z(G)\Lambda$, 
and $TZ(G)\Lambda$ is a closed subgroup of $G$ containing the lattice 
$\Lambda$.  Thus $TZ(G)\Lambda/\Lambda$ has finite volume, which in turn
implies that $TZ(G)\cap\Lambda$ is a lattice in $TZ(G)$. Since 
$Z(G)$ is infinite by hypothesis, and since $T$ is a torus, we conclude that
$\Lambda' := TZ(G)\cap\Lambda$ is an infinite normal abelian subgroup
of $\Lambda$. Since $[\Lambda, \Lambda'] \subset [\Lambda, T]\subset T$,
and since $T$ is compact, we have that $[\Lambda, \Lambda']$ 
is finite. Now since $\Lambda$ is
finitely generated (every lattice in a connected Lie group is
finitely generated), we can conclude easily that a subgroup of 
$\Lambda'$ of finite index is contained in $Z(\Lambda)$. This
proves that $Z(\Lambda)$ is infinite.
\endproof

\section{Some applications}
\label{section:consequences}

In this section we finish the proof of Theorem \ref{theorem:locsym1}.
We then use Theorem \ref{theorem:main} and its proof, and also 
Theorem \ref{theorem:locsym1}, to prove the other theorems and
corollaries stated in the introduction.

\subsection{No normal abelian subgroups (proof of Theorem \ref{theorem:locsym1})}
\label{section:proof:locsym}

The fact that (2) implies (1) follows immediately from well-known
properties of closed, locally symmetric Riemannian manifolds.  Such $M$
are aspherical by the Cartan-Hadamard theorem.  Any normal abelian
subgroup is trivial since the symmetric space $\widetilde{M}$ has no
Euclidean factors.  The other two properties follow from the
definitions.

To prove that (1) implies (2), we first quote Proposition 
\ref{prop:no:normal:abelian} followed by Proposition 
\ref{proof:semisimp:consequence}. This gives that $M$ has a 
finite-sheeted Riemannian cover $M'$ of $M$ which is a smooth (indeed 
Riemannian warped) product $M'=N\times B$, where $N$ is 
is isometric to a nonempty, irreducible, locally
symmetric, nonpositively curved manifold.  
But $M'$ is smoothly irreducible by hypothesis, so
that $B$ must be a single point.  It follows that $M'=N$ is locally
symmetric.  Since the metric on $M'$ was lifted from $M$, we have that
$M$ is locally symmetric.

\subsection{Word-hyperbolic groups 
(proof of Corollary \ref{corollary:locsym3})}

Again, (2) implies (1) follows immediately from the basic properties of
closed, rank one locally symmetric manifolds.

To prove that (1) implies (2), first note that no torsion-free
word-hyperbolic group can virtually be a nontrivial product, since then
it would contain a copy of $\Z\times \Z$.  It then follows from Theorem
\ref{theorem:locsym1} that $M$ is locally symmetric.  
But every closed, locally symmetric manifold $M$ either 
contains $\Z\times \Z$ in its fundamental group, or $M$ must be negatively
curved; hence the latter must hold for $M$.

\subsection{Almost simple groups (proof of Corollary \ref{corollary:locsym2})}

This follows just as the proof of Corollary \ref{corollary:locsym3}, 
but using the following fact: an irreducible, 
cocompact lattice in a noncompact 
semisimple Lie group $G$ is almost simple if and
only if $\rank_{\R}\geq 2$.  The ``if'' direction is the statement of
the Margulis Normal Subgroup Theorem (see \cite{Ma}, Thm. IX.5.4).  
For the ``only if'' direction, first recall that cocompact lattices
in rank one semisimple Lie groups are non-elementary word-hyperbolic.
Such groups are never almost simple; for example, a theorem 
of Gromov-Olshanskii (see \cite{Ol}) gives that all such groups 
have infinite torsion quotients. 

\subsection{Universal bound (proof of Theorem \ref{theorem:magic:number})}

A theorem of Kazhdan-Margulis (see, e.g., \cite{Ra}, Corollary XI.11.9) 
shows that, for every connected 
semisimple Lie group $G$, 
there exists $\epsilon=\epsilon(G)$ such that the
covolume of every lattice in $G$ is greater than $\epsilon$.  Let 
$\epsilon$ be this constant for $G=\Isom(\widetilde{M})$.  Now let 
$C_1=\vol(M,g_0)/\epsilon$.  

Now let $X$ denote $\widetilde{M}$ endowed with any fixed Riemannian
metric $h$ lifted from $M$.  If $\Isom(X)$ is not discrete, then by
Theorem \ref{theorem:locsym1} we have that $h$ is homothetic to $g_0$,
and we are done.  If $\Isom(X)$ is discrete, then $X/\Isom(X)$ is a 
compact orbifold, Riemannian covered by the compact manifold $M$, with
degree of the cover $d:=[\Isom(X):\pi_1(M)]$.  Since the volume of a
cover is multiplicative in degree, we have that 
$$\vol(M,h)=d\vol(X/\Isom(X)))>d\epsilon$$

Now, Gromov has shown (see \cite{Gr}, and also \cite{CF} for a
different, and more detailed, proof) that the ``minvol'' 
invariant of $M$ is positive; this means that
there is a constant $C_2=C_2(M)$ so that $\vol(M,h)\leq C_2\vol(M,g_0)$
for all $h$ with sectional curvatures bounded above by $1$ in absolute
value.  After rescaling, we can assume the given $h$ satisfies this
curvature bound.  We then have 
$$d<\vol(M,h)/\epsilon\leq C_2\vol(M,g_0)/\epsilon=C_1C_2$$
and we are done after setting $C=C_1C_2$.

\subsection{Models for compact and finite volume manifolds 
(proof of Theorem \ref{theorem:covering})}
\label{section:lattices}

Let $\Gamma_1$ (resp. $\Gamma_2$) be a cocompact (resp. noncocompact)  
lattice in $\Isom(X)$.  Since $\Isom(X)$ contains $\Gamma_1$, it acts 
cocompactly on $X$.  Suppose $\Isom(X)$ were discrete, so that
$\Isom(X)$ acts properly and cocompactly on $X$.  Then the 
(orbifold) quotient $\Isom(X)\backslash X$ would have finite volume.  Since
covolume is multiplicative in index, it would follow that
$\Gamma_2<\Isom(X)$ has finite index.  But then $\Gamma_2$ would act
cocompactly since $\Isom(X)$ does, a contradiction.  Hence $\Isom(X)$ is
not discrete and we may apply Theorem \ref{theorem:main}.

We thus obtain a Riemannian orbibundle, which at the level of
universal covers gives a Riemannian warped product structure
$Y\tilde{\times} X_0$, with $Y$ the universal cover of $B$, where the
metric has the property that for each $x\in X$, the metric on $x\times
X_0$ is an $I_0$- homogeneous metric, depending on $x$.

Let $\Lambda_i:=\Gamma_i\cap I_0, i=1,2$.   Claim I of the proof of Theorem
\ref{theorem:main} gives that $\Lambda_1$ is cocompact.  We must now
prove that $\Lambda_2$ is a noncocompact lattice.

To this end, first note that $\pi:=\Gamma_2/\Lambda_2$ is a cocompact 
lattice in $\Isom(Y)$.  First suppose $B$ is not $1$-dimensional.  
We can then perturb the
metric on $B$ to get a new universal (in the category of orbifolds) 
cover $Y'$ with
$\Isom(X)=\Isom(Y'\times X_0)$ but with $\Isom(Y')=\pi_1^{orb}(B)$.  Now
$\pi$ is a lattice in $\Isom(Y')$, so it has finite index in
$\pi_1(B)$.  

We thus have that each of $\Gamma_i, i=1,2$ 
can be written as a group extension with 
kernel $\Lambda_i$ and quotient a group with the same rational 
cohomological dimension $\qcd$ as $\pi_1(B)$.  We consider 
rational cohomological dimension $\qcd$ in order to deal with the fact that 
$\pi_1(B)$ might not be virtually torsion-free.

Since each $\Gamma_i$ 
acts properly on the contractible manifold $X$, and since 
$\Gamma_1$ acts cocompactly and $\Gamma_2$ does not, we have that 
$\qcd(\Gamma_2)<\qcd(\Gamma_1)$ (see, e.g. \cite{Bro}, VIII.8).  

Now, if $\Lambda_2$ were cocompact,
then $\Gamma_2$ would be an extension of fundamental groups of 
closed, aspherical manifolds, and so $\qcd(\Gamma_2)$ would be the sum
of the $\qcd$ of the kernel and quotient; but this sum equals
$\qcd(\Gamma_1)$ (see, e.g.,  Theorem 5.5 of \cite{Bi}), 
a contradiction.  Hence $\Lambda_2$ is not
cocompact.  

If $B$ is one-dimensional, then no perturbation as above exists.  To
remedy this, we simply take the product of $B$ with a closed, 
genus $2$ surface, endowed with a Riemannian metric with trivial
isometry group.  We then run the rest of the argument verbatim.

\subsection{Irreducible lattices (proof of Theorem \ref{theorem:irreducible})}
\label{section:irreducible}

We first note that the hypothesis implies $\Isom(X)$ is not discrete.  
Since in addition $X$ has a cocompact discrete subgroup, we may apply Theorem
\ref{theorem:main}.  Hence $X$ is isometric to a warped Riemannian
product $X=Y\tilde{\times}X_0$. By the proof of Theorem
\ref{theorem:main}, the group $\Isom(X_0)$ corresponds with connected
component of the identity of $\Isom(X)$.  In particular $\Isom(Y)$ must
be discrete.  The theorem follows easily.

\subsection{The Hopf Conjecture}
\label{section:Hopf}

A well-known conjecture of Hopf-Chern-Thurston
states that the Euler characteristic of any 
closed, aspherical manifold $M^{2k}$ 
satisfies $(-1)^k\chi(M^{2k})\geq 0$.  A stronger conjecture of Singer
posits that the $L^2$-cohomology of $M^{2k}$ vanishes except in dimension
$k$ (see, e.g. \cite{CG,Lu}).  These conjectures are completely 
open except when $k=1$.

We will now prove that the Hopf Conjecture holds for those smooth,
aspherical manifolds $M^{2k}$ which admit {\em some} Riemannian metric with
symmetry. 

\begin{theorem}
\label{theorem:Hopf}
Let $M^{2k}$ be any closed, aspherical, smooth manifold which is 
smoothly irreducible.  If $M$ admits some Riemannian metric so
that the induced metric on the universal cover $\widetilde{M}$ satisfies
$[\Isom(\widetilde{M}):\pi_1(M)]=\infty$, 
then the Singer Conjecture (and hence
the Hopf Conjecture) is true for $M^{2k}$.
\end{theorem}

The Singer Conjecture clearly holds for products 
of surfaces and also for products of any $3$-manifolds with $S^1$.  Thus
Theorem \ref{theorem:Hopf} holds in dimension four without
the assumption that $M^4$ is smoothly irreducible.

\bigskip
\noindent
{\bf Proof of Theorem \ref{theorem:Hopf}. }
If $\pi_1(M^{2k})$ has no nontrivial normal abelian subgroups and is smoothly
irreducible, then by Theorem
\ref{theorem:locsym1} it admits a locally symmetric Riemannian 
metric.  The Singer Conjecture is known for such manifolds (see,
e.g. \cite{Lu}, Cor. 5.16).  If $\pi_1(M^{2k})$ does contain a nontrivial,
normal abelian subgroup $A$, then one may apply a theorem 
of Cheeger-Gromov (\cite{CG}, Cor. 0.6) 
which gives, even more generally for amenable 
$A$, that Singer's Conjecture holds.
\endproof

\subsection{Complex manifolds (proof of Theorem \ref{theorem:frankel})}
\label{section:complex}

As pointed out above, Theorem \ref{theorem:frankel} follows
immediately from Proposition 0.1 in \cite{Na} together with 
Claim IV in the proof of Theorem \ref{theorem:main}.  

Actually, we can do without Nadel's result under certain mild
assumptions.  Since $\Aut(\widetilde{M})$ acts isometrically in some
Riemannian metric, we can apply Theorem \ref{theorem:main} directly to
obtain the claimed splitting, but with two problems: first, the
splitting is an isometric (not holomorphic) one; and second, the factor
$M_1$ is only locally homogeneous, not necessarily locally symmetric.  

The second problem is corrected once we know that $\Aut(\widetilde{M})^0$ 
is semisimple with finite center; equivalently, $M_1$ doesn't fiber with
nontrivial solvmanifold fiber.  One way to rule this out is to assume
that $\pi_1(M)$ contains no infinite normal abelian subgroups, giving us
a new proof of Theorem \ref{theorem:frankel} in this case.  In complex
dimension two, $M_1$ must be locally symmetric, for otherwise it would 
have a fibering with torus fiber, giving that $\chi(M)=0$.  However, in
complex dimension two, we have 
$$\chi(M)=c_2(M)\geq \frac{1}{3}c_1(M)^2>0$$
by the Bogomolov-Miyaoka-Yau inequality.

To correct the first problem, once we have that the factor $M_1$ is
locally symmetric, we apply Siu's Rigidity Theorem (also used in
\cite{Fr2}) to obtain that $M_1$ is biholomorphic to a closed, 
Hermitian locally symmetric space.

\section{An extension to the non-aspherical case}
\label{section:nonaspherical}

A number of the results from this paper generalize from closed, aspherical 
Riemannian manifolds to all closed Riemannian manifolds. 
As an illustrative example we give 
the following theorem.  As far as we know, this is the
first geometric rigidity theorem for non-aspherical manifolds with
infinite fundamental group.

We say that the universal cover $\widetilde{M}$ of a Riemannian manifold
$M$ has {\em essential extra symmetry} if, for a compact subset
$K\subset\widetilde{M}$, and for 
all $\epsilon >0$, there exists $g\in\Isom(\widetilde{M})$ 
such that
\begin{itemize}
\item $d(g(m),m)<\epsilon$ for all $m\in K$
\item $\sup d(g^n(m),m)=\infty$ for any fixed $m$.
\end{itemize}

The condition that $\widetilde{M}$ have an essential symmetry is 
equivalent to the identity component of $\Isom(\widetilde{M})$ being
noncompact.  If $\pi_1(M)$ is torsion free, this is equivalent to
$\Isom(\widetilde{M})$ not being an extension of $\pi_1(M)$ by any compact
group.  (If the identity component of $\Isom(\widetilde{M})$ were
compact, then $\pi_1(M)$ would
intersect this group in a normal lattice by Claim I, and this 
must be trivial if $\pi_1(M)$ is torsion free.)

\begin{theorem}
\label{theorem:nonaspherical}
Let $M$ be any closed Riemannian manifold which, for simplicity, is 
smoothly irreducible.  Suppose that $\pi_1(M)$
contains no infinite, finitely-generated, 
normal abelian subgroup, and that $M$ has
essential extra symmetries. Then there exists a finite Riemannian 
cover $M'$ of $M$ which is a fiber bundle over a closed, irreducible,
locally symmetric manifold, with all fibers isometric.  In particular
the structure group of the bundle is compact.
\end{theorem}

The proof of Theorem \ref{theorem:nonaspherical} follows 
quite closely the proof of Theorem \ref{theorem:locsym1}, with only a
few adaptations.  Hence instead of a detailed proof, we now just indicate the
adaptations that are necessary.

\bigskip
\noindent
{\bf Proof. }We will describe how to modify the proof of 
Theorem \ref{theorem:locsym1} in order to prove Theorem
\ref{theorem:nonaspherical}.  As in the proof of Theorem
\ref{theorem:locsym1}, the isometry group $\Isom(\widetilde{M})$ is a semidirect product of $I_0$ and $\Delta$, where
$\Gamma=\pi_1(M)$ intersects $I_0$ in a
lattice $\Gamma_0$.  As before, the hypothesis that $\Gamma$ has no
infinite, finitely-generated, normal abelian subgroup implies that $I_0$
is semisimple with finite center.

This implies as in \S\ref{section:when:semisimp} that, replacing $M$ by
a finite cover if necessary, we have a product structure on $\Gamma$.
We then have a harmonic map $f:M \to K_0\backslash I_0/\Gamma_0$, where
$K_0$ is a maximal compact subgroup of $I_0$.  
This will provide the orbibundle structure on $M$, and the Lie group
$I_0$ will be responsible for the
isometries among the fibers.

The proof of this is again based on the uniqueness of harmonic maps to
nonpositively curved manifolds, when one has a surjection of
fundamental groups, and no center in the image.  Here, since $\Gamma$
splits as a product, having such a center would
give us a normal abelian subgroup, contradicting the hypothesis.  
On $\widetilde{M}$ the $\Gamma$-equivariant map given by the lift of
$f$ extends, by the arithmetic
group trick of Claim \ref{claim:equiv}, to an $I_0$-equivariant 
map $F:\widetilde{M}\to I_0/K_0$.  
Consequently the fibers are isometric to each other as the 
$I_0$-action on the target is transitive. 
\endproof

\section{A truly singular orbibundle}
\label{section:example:nocover}

In this section we construct a ($7$-dimensional) 
Riemmanian manifold $M$ with the property
that $M$ is a Riemannian orbibundle, but no finite cover of $M$ is a
Riemannian fiber bundle.  $M$ will also have the property that
$\Isom(\widetilde{M})$ is not discrete.  This will prove that the
``orbibundle'' conclusion of Theorem
\ref{theorem:main} cannot be improved to a fiber bundle, even if one is
willing to pass to a finite cover.
  
We begin with the following construction, which we believe is of
independent interest.

\begin{theorem}
\label{theorem:weirdaction}
There is a group $\Gamma$ which acts cocompactly, properly
discontinuously by diffeomorphisms on $\R^n$ for some $n$, 
but which is not virtually
torsion free. Moreover,  there exists $\xi\in H^2(\Gamma,\Z)$ 
which restricts to a nonzero class on some $\Z/p\Z \in \Gamma$, where $p$ is a
prime, and which also restricts to a nonzero class on every finite index
subgroup. 
\end{theorem}     

\proof     
Let $G= \Z/p\Z\ast \Z$.  Then $G$ acts properly discontinuously and
cocompactly on its Bass-Serre tree $T$.   In \cite{Bri}, Bridson 
shows that there is an amalgamated free
product $H=G\ast_FG$, with $F$ a nonabelian free group, 
which has no finite quotients.  Let $Y$ be the universal cover
of the standard Cayley
$2$-complex for $H$; hence $H$ acts properly discontinuously and
cocompactly on $Y$.  It is easy to see that $Y$ is contractible. 

Equivariantly thicken $Y$ (see \cite{As} for the details of equivariant 
thickening).  The action on this thickening 
gives a properly discontinuous, cocompact 
$\Gamma$-action on a contractbile $6$-dimensional manifold.  Since the
action is simplicial, we can perform the Davis reflection group
construction (see \cite{Da}) equivariantly to build a cocompact action
of a group $\Lambda$ on a contractible manifold, and with $H$ being a
retract of $\Lambda$.  This action can in
fact be made smooth, as is explained in
\S 17 of \cite{Da}.  By a theorem of Stallings \cite{St}, 
the group $\Gamma:=\Lambda\times \Z$ then acts on $\R^n$ properly
discontinuously and cocompactly by diffeomorphisms.

Note that any finite index subgroup of $\Gamma$  intersects $H$ 
in a finite index subgroup of $H$, which must therefore 
be all of $H$ since $H$ has no finite index subgroups.  
Hence the fixed $\Z/p\Z$ subgroup of $\Gamma$ must lie in each finite
index subgroup of $\Gamma$.  

Now by construction there is a surjection $\Gamma\to G$.  Let $\xi\in
H^2(\Gamma,\Z)$ denote the pullback of the class generating
$H^2(G,\Z)$.  As the amalgamating subgroup $F$ is free, the amalgamation
dimension of $H^2(H,\Z)$ is at most one greater than 
$H^2(G,\Z)$.  Since $H$ is a retract of $\Lambda$, we know that 
$\xi \in H^2(\Lambda,\Z)$ is also nonzero.  It follows easily that 
$\xi$ pulls back to a nonzero class in $\Gamma$.
\endproof

We now build the manifold $M$.  Let $\Gamma, \xi$ be given as in 
Theorem \ref{theorem:weirdaction}.  Let 
$\widetilde{\Gamma}$ denote the central extension of $\Gamma$ given
by the cocycle $\xi\in H^2(\Gamma,\Z)$.  Since this cocycle vanishes in 
$H^2(\Gamma,\R)$, we have that $\widetilde{\Gamma}$ lies in 
$\Gamma \times\R$.  Now fix any $\Gamma$-invariant metric on $\R^n$,
and extend this to any $(\Gamma\times \R)$-invariant metric on
$\R^n\times \R=\R^{n+1}$; call the resulting Riemannian manifold $Y$.  
 
Now $\Gamma\times \R$ acts properly discontinuously and 
cocompactly by isometries on $Y\approx \R^{n+1}$.  The quotient $M$
clearly satisfies the claimed properties.  Note too that
$\Isom(\widetilde{M})=\Isom(Y)$ contains $\R$, and so is not discrete.

\medskip
\noindent
{\bf Remark. }in the examples above, the dimension of $M$ is at least
$7$.  We
do not know whether this dimension can be lowered.  Indeed, it seems
difficult to obtain information about the geometry of such examples,
although they do seem compatible with at least large-scale
nonpositive curvature.

\bigskip
\noindent
Dept. of Mathematics, University of Chicago\\
5734 University Ave.\\
Chicago, Il 60637\\
E-mail: farb@math.uchicago.edu, shmuel@math.uchicago.edu

\end{document}